\documentclass{article}

\usepackage{amsmath}
\usepackage{graphicx}
\usepackage{amsfonts}
\usepackage{amssymb}
\usepackage{psfrag}
\usepackage{amscd}

\newcommand{\re}{\mathbb{R}}
\newcommand{\cpx}{\mathbb{C}}

\newcommand{\N}{\mathbb{N}}

\newcommand{\lmd}{\lambda}

\newcommand{\eps}{\epsilon}

\newcommand{\sig}{\sigma}

\newcommand{\proof}{\noindent {\it Proof. \quad  }}
\newcommand{\eproof}{{\ $\square$ }}
\newcommand{\reff}[1]{(\ref{#1})}
\newcommand{\pt}{\partial}

\newcommand{\bdes}{\begin{description}}
\newcommand{\edes}{\end{description}}

\newcommand{\bal}{\begin{align}}
\newcommand{\eal}{\end{align}}

\newcommand{\bnum}{\begin{enumerate}}
\newcommand{\enum}{\end{enumerate}}

\newcommand{\bit}{\begin{itemize}}
\newcommand{\eit}{\end{itemize}}

\newcommand{\bea}{\begin{eqnarray}}
\newcommand{\eea}{\end{eqnarray}}
\newcommand{\be}{\begin{equation}}
\newcommand{\ee}{\end{equation}}

\newcommand{\baray}{\begin{array}}
\newcommand{\earay}{\end{array}}

\newcommand{\bsry}{\begin{subarray}}
\newcommand{\esry}{\end{subarray}}

\newcommand{\bca}{\begin{cases}}
\newcommand{\eca}{\end{cases}}

\newcommand{\bcen}{\begin{center}}
\newcommand{\ecen}{\end{center}}

\newcommand{\bbm}{\begin{bmatrix}}
\newcommand{\ebm}{\end{bmatrix}}

\newcommand{\bmx}{\begin{matrix}}
\newcommand{\emx}{\end{matrix}}

\newcommand{\bpm}{\begin{pmatrix}}
\newcommand{\epm}{\end{pmatrix}}

\newcommand{\btab}{\begin{tabular}}
\newcommand{\etab}{\end{tabular}}

\newsavebox{\myname}
\newsavebox{\emailaddr}

\sbox{\myname}{\bf Jiawang Nie}   
\sbox{\emailaddr}{email: njw@math.berkeley.edu}         

\newcommand{\thmlist}{
\begin{list}{Step 1}
{\setlength{\leftmargin}{0.6 in}\setlength{\labelwidth}{0.5 in}} }
\newcommand{\alglist}{
\begin{list}{Step 1}
{\setlength{\leftmargin}{1.1 in}\setlength{\labelwidth}{1.0 in}} }
\newtheorem{thm}{Theorem}[section]

\newtheorem{lemma}[thm]{Lemma}

\newtheorem{exm}[thm]{Example}
\newcommand{\Section}[1]{
    \par
    \stepcounter{section}
    \settowidth{\hangindent}{\large\bf\thesection.~}
    \hangafter=1
    \bigskip\bigskip\noindent
    {\large\bf\hbox{\thesection.~}#1}\par
    \nopagebreak
    \medskip
    \renewcommand{\theequation}{\thesection.\arabic{equation}}
    \setcounter{equation}{0}
    \setcounter{subsection}{0}
}
\renewcommand{\subsection}[1]{
    \stepcounter{subsection}
    \settowidth{\hangindent}{\bf\thesubsection.~}
    \hangafter=1
    \bigskip\bigskip\noindent
    {\bf\hbox{\thesubsection.~}#1}\par
    \nobreak
    \medskip
}
\usepackage{amsmath}
\usepackage{graphicx}
\usepackage{colortbl}
\usepackage{color}
\usepackage{amsfonts}
\usepackage{amssymb}
\usepackage{psfrag}
\usepackage{amscd}

\begin{document}

\title{ Minimizing Polynomials Over Semialgebraic Sets
\footnote{This paper was supported in part by National Science Foundation (ELA-0122599)}
\author{ Jiawang Nie\footnote{{\tt njw@math.berkeley.edu}, Department of Mathematics,
University of California, Berkeley, CA 94720},
\ James W. Demmel\footnote{{\tt demmel@cs.berkeley.edu}, Department of Mathematics
and EECS, University of California, Berkeley, CA 94720},
and\ Victoria Powers\footnote{{\tt vicki@mathcs.emory.edu},
Deparment of Mathematics and Computer Science, Emory University,
Atlanta, GA 30322}
\date{}
}}

\maketitle

\abstract{
This paper concerns a method for finding the minimum of a
polynomial on a semialgebraic set, i.e.,  a set in $\re^m$ defined
by finitely many polynomial equations and inequalities,
using the Karush-Kuhn-Tucker (KKT) system
and  sum of squares (SOS) relaxations.
This generalizes results in the recent paper \cite{njw_grad},
which considers  minimizing polynomials on
algebraic sets, i.e.,  sets in $\re^m$ defined by finitely many
polynomial equations.
Most of the theorems and conclusions in \cite{njw_grad} generalize
to semialgebraic sets, even in the case where
the semialgebraic set is not compact.
We discuss the  method in some special cases, namely,
when the semialgebraic set is contained in
the nonnegative orthant $\re^n_+$ or in box constraints $[a,b]_n$.
These constraints make the computations more efficient.
}

\hspace{2cm}


\Section{Introduction}

In this paper, we consider the optimization problem
\begin{align}
f^*=\min &\,\,\,\ f(x)  \label{opt1} \\
s.t. &\,\,\,\ g_i(x)=0,\,\,\ i=1,\cdots,s, \label{opt2}\\
     &\,\,\,\ h_j(x)\geq 0,\,\,\  j=1,\cdots,t \label{opt3}
\end{align}
where $x=\bbm x_1 & \cdots x_n\ebm^T \in \re^n$
and $f(x),g_i(x),h_j(x)\in \re[x]$
(the ring of real multivariate polynomials in $x$).
Let ${\cal F}$ be the feasible region, i.e., the subset of $\re^n$
which satisfies constraints $\reff{opt2}-\reff{opt3}$; ${\cal F}$
is a semialgebraic set.
Many optimization problems in practice
can be formulated as \reff{opt1}-\reff{opt3}.
Finding the global optimal solutions to
$\reff{opt1}-\reff{opt3}$ is an NP-hard problem,
even if $f(x)$ is quadratic and $g_i,h_j$ are linear.
For instance, the Maximum-Cut problem for graphs is of this form,
and it is NP-hard (\cite{NP_complete}).

Recently, the techniques of sum of squares (SOS) relaxations
and moment matrix methods
have made it possible to find the global optimal solutions
to \reff{opt1}-\reff{opt3}
by approximating nonnegative polynomials
with SOS polynomials,
which allows the problem to be implemented as a semidefinite program.
For more details about these methods and their
applications, see
\cite{lasserre, Laurent, njw_elip, njw_grad, Pths, Par00, Par01,Par02}.
To prove the convergence of these methods,
it is often necessary to  assume that the feasible region ${\cal F}$
is compact or even finite.
In \cite{Par02}, it is shown that SOS relaxations can
solve \reff{opt1}-\reff{opt3} globally in finitely many steps
in the case where
$\{x\in\cpx^n:\,\ g_1(x)=\cdots=g_s(x)=0\}$
is finite and the ideal
$\langle g_1(x),\cdots,g_s(x)\rangle$
is radical.
If we only assume that
$\{x\in\cpx^n:\,\ g_1(x)=\cdots=g_s(x)=0\}$
is finite, it is shown in
\cite{Laurent}  that the moment matrix method can
solve \reff{opt1}-\reff{opt3} globally
in finitely many steps.  Finally, if
${\cal F}$ is compact and the set of polynomials $\{ g_i, h_i \}$
satisfies an additional assumption (see Theorem~\ref{putinar}),
then arbitrarily close lower bounds for $f^*$
can be obtained by SOS relaxations or moment matrix methods \cite{lasserre}.
In this case,  the convergence is asymptotic, however
little is known about the errors in the bounds.

The above global optimization methods are based on
representation theorems from real algebraic geometry for polynomials
positive and nonnegative on semialgebraic sets.
On the other hand, the traditional local methods
in optimization often follow the
first order optimality conditions
(zero gradient in the unconstrained case
or the Karush-Kuhn-Tucker (KKT) system
in the constrained case).
The underlying idea in \cite{njw_grad} and the present paper is
to combine these two types of methods in order to more efficiently
solve
\reff{opt1}-\reff{opt3} globally.
In \cite{njw_grad},
SOS relaxations are applied
on the gradient ideal ${\cal I}_{grad}$
(the ideal generated by all the partial
derivatives of $f(x)$) in the unconstrained case,
and on the KKT ideal $I_{KKT}$ in the constrained case, where
only equality constraints are allowed.
When ${\cal I}_{grad}$ or $I_{KKT}$ is radical,
which is generically true in practice,
the method in \cite{njw_grad}
can solve the optimization \reff{opt1}-\reff{opt2} globally;
otherwise, arbitrarily close lower bounds
of $f^*$ can be obtained.
No assumptions about
${\cal F}$ are made, i.e., it need not be finite or even compact.

The convergence of the method in \cite{njw_grad}
assumes that the constraints are algebraic sets.
If there are any inequality constraints,
${\cal F}$ is no longer algebraic
but only semialgebraic and the proof in \cite{njw_grad} does not work.
The motivation of this paper is to
generalize the method in \cite{njw_grad}
to handle semialgebraic constraints.

The KKT system of problem \reff{opt1}-\reff{opt3} is
\begin{align}
F \overset{\Delta}{=}
\nabla f(x) &+ \sum_{i=1}^s \lmd_i \nabla g_i(x)
-\sum_{j=1}^t \nu_j \nabla h_j(x) =0, \label{kkt1}\\
h_j(x)&\geq 0,\nu_j h_j(x)=0,\,\ j=1,\cdots,t, \label{kkt2}\\
g_i(x)& = 0,\,\, i=1,\cdots,s, \label{kkt3}
\end{align}
where  vectors $\lmd=\bbm \lmd_1 & \cdots \lmd_s\ebm^T$ and
$\nu=\bbm \nu_1 & \cdots \nu_t\ebm^T$
are called Lagrange multipliers.
See \cite{Nocedal} for some regularity conditions
that make the KKT system hold
at local or global minimizers.
For an example where the KKT system fails,
see Example~\ref{kktfails} in Section~4.

Note that we do not require $\nu\geq 0$ in the above;
this makes the SOS relaxations simpler and
does not affect the convergence of the method, since
omitting the constraint $\nu \geq 0$ means simply that
there are more feasible points for \reff{kkt1}-\reff{kkt3},
including maxima as well as minima.
But since we minimize over this larger set, we get
the same minima. Minimizing over this larger set
makes our problem easier, because it reduces the number
of inequality constraints, which as we will see greatly
lowers the complexity of our algorithm.

Let $f^*_{KKT}$ be the global minimum of
$f(x)$ over the KKT system defined by \reff{kkt1}-\reff{kkt3}.
Assume the KKT system holds at the global minimizers.
Then we claim that $\mathbf{f^*=f^*_{KKT}}$.
First, $f^*\leq f^*_{KKT}$ follows immediately
from the fact that all $x$ in the KKT system are feasible.
Now let $x^*$ be a global minimizer such that
$f(x^*)=f^*$, then by assumption, there exist Lagrange multipliers
$\lmd^*$ and $\nu^*\geq 0$ such that
$(x^*,\lmd^*,\nu^*)$ satisfies the above KKT system.
Thus $f^*\geq f^*_{KKT}$ and hence they are equal.

Define the KKT ideal $I_{KKT}$ and its varieties as follows:
\begin{align*}
I_{KKT}&=\left< F_1,\cdots,F_n, g_1,\cdots,g_s,\nu_1h_1,\cdots,\nu_t h_t\right>, \\
V_{KKT}&=\{(x,\lmd,\nu)\in \cpx^n\times \cpx^s\times \cpx^t:
p(x,\lmd,\nu)=0,\,\,\forall p \in I_{KKT}\}, \\
V_{KKT}^\re&=\{(x,\lmd,\nu)\in \re^n\times \re^s\times \re^t:
p(x,\lmd,\nu)=0,\,\,\forall p \in I_{KKT}\}.
\end{align*}
Here $F=[F_1,\cdots,F_n]^T$ is defined in \reff{kkt1}.
Let $${\cal H}=\{(x,\lmd,\nu)\in \re^n\times \re^s\times \re^t:
h_j(x)\geq 0,\,\ j=1,\cdots,t \}.$$
The {\it preorder} cone $P_{KKT}$ associated with the KKT system is
$$
P_{KKT}=
\left\{ \left. \sum_{\theta \in \{0,1\}^t}
\sig_{\theta} h_1^{\theta_1}  h_2^{\theta_2} \cdots h_t^{\theta_t} \right| \,\
\sig_{\theta} \mbox{ are SOS}\,\right\}
+I_{KKT}.
$$
The {\it linear} cone associated with the KKT system is
$$
M_{KKT}=
\left\{ \left. \sig_0
+\sum_{j=1}^t \sig_jh_j \right| \,\
\sig_0,\cdots,\sig_{t} \mbox{ are SOS}\,\right\}
+I_{KKT}.
$$
Note that $I_{KKT}\subseteq
M_{KKT} \subseteq
P_{KKT }\subseteq \re[x,\lmd,\nu]$.

In solving SOS programs,
we often set an upper bound on the degrees of the involved polynomials.
Define the truncated KKT ideal
\begin{align*}
I_{N,KKT}=\Big\{
\sum_{k=1}^n  \phi_k F_k+
&\sum_{i=1}^s \varphi_i g_i+\sum_{j=1}^t  \psi_j \nu_j h_j\Big|
deg(\phi_k F_k), deg(\varphi_i g_i), deg(\psi_j\nu_j h_j) \leq N
\Big\}.
\end{align*}
and truncated preorder and linear cones
\begin{align*}
P_{N,KKT}& = \left\{
\left.\sum_{\theta \in \{0,1\}^t}
\sig_{\theta} h_1^{\theta_1}  h_2^{\theta_2} \cdots h_t^{\theta_t} \right|
deg(\sig_{\theta} h_1^{\theta_1}\cdots h_t^{\theta_t}) \leq N
\right\}
+I_{N,KKT}.  \\
M_{N,KKT}&= \left\{ \sig_0
+\sum_{j=1}^t \sig_jh_j
\left|
\bmx
\sig_0,\cdots,\sig_{t} \mbox{ are SOS} \\
deg(\sig_0), deg(\sig_jh_j) \leq N \\
\emx \right.
\right\}
+I_{N,KKT}.
\end{align*}

A sequence $\{p^*_N\}$ of lower bounds of \reff{opt1}-\reff{opt3}
can be obtained by the following SOS relaxations:
\begin{align}
p^*_N= \max_{\gamma \in \re} &\,\,\,\ \gamma  \label{sospod1} \\
s.t. &\,\,\,\ f(x)-\gamma \in P_{N,KKT}. \label{sospod2}
\end{align}
Since $P_{N,KKT}$ has a summation over $2^t$ terms like
$\sig_{\theta} h_1^{\theta_1}  h_2^{\theta_2} \cdots h_t^{\theta_t}$,
it is usually very expensive to solve the SOS program
\reff{sospod1}-\reff{sospod2} in practice.
So in practice, it is natural to replace the truncated
preorder cone $P_{N,KKT}$ by truncated linear cone  $M_{N,KKT}$,
which leads to the SOS relaxations:
\begin{align}
f^*_N= \max_{\gamma \in \re} &\,\,\,\ \gamma  \label{sosrlx1} \\
s.t. &\,\,\,\ f(x)-\gamma \in M_{N,KKT}. \label{sosrlx2}
\end{align}
Thus we have the increasing sequences of lower bounds
$\{f^*_N\}_{N=2}^\infty$ and $\{p^*_N\}_{N=2}^\infty$ such that
$f^*_N \leq p^*_N \leq f^*.$

The following notation is used throughout:
We denote by $deg(p)$ the degree of a polynomial $p$.
%
%
The vector inequality $u\leq v\,(u,v\in \re^n)$
is defined component-wise, i.e.,
$u_i\leq v_i$ for each $i$.
$[u,v]_n$ denotes the set of all vectors $w\in \re^n$
such that $u\leq w \leq v$.

This paper is organized as follows.
Section~2 is a review of some fundamental results
from algebraic geometry.  In
Section~3 we discuss the representation of the
polynomial $f(x)$ in the cones  $M_{KKT}$ and $P_{KKT}$.
We analyze the convergence
of the lower bounds $\{p^*_N\}$ and $\{f^*_N\}$ in Section~4.
In Section~5, we consider
some special cases of inequality constraints, in particular,
the nonnegative orthant $\re^n_+$ and the box $[a,b]_n$.
Section~6 draws conclusions.

\Section{Preliminaries}

This section will introduce some basic notions from
algebraic geometry needed for our discussion.
Readers may consult \cite{CLO97,CLO98,Eisenbud}
for more details.
In this section, all polynomials are
in the indeterminate $x=(x_1,\ldots,x_m)$
for the simplicity of notation.
Here $x$ is {\em not} the ``$x$'' in the Introduction,
but rather a generic indeterminate.
In later sections, $x$ will be again
the ``$x$'' in \reff{opt1}-\reff{opt3},
and all polynomials will be in the variables $(x,\lmd,\nu)$,
unless explicitly stated otherwise.

We write $\,\re[x] = \re[x_1,\ldots,x_m]$
for the ring of polynomials in indeterminates
$x = (x_1,\ldots,x_m)$ with real coefficients.
A polynomial $p \in \re[x]$ is {\it SOS} if it can be
written as a sum of squares of polynomials in $\re[x]$.
A subset $I$ of $\re[x]$ is an {\it ideal}
if $p\cdot q\in I$ for any $p\in I$ and $q\in \re[x]$.
For $p_1,\ldots,p_r \in \re[x]$,
$ \langle p_1,\cdots,p_r \rangle $ denotes
the smallest ideal containing the $p_i$. Equivalently,
$ \langle p_1,\cdots,p_r \rangle $ is the set of all polynomials that are
polynomial linear combinations of the $p_i$.
Every ideal arises in this way:

\begin{thm}[Hilbert Basis Theorem]  \label{HilBasis}
\hfill \break
Every ideal $I\subset \re[x]$
has a finite generating set, i.e.,
$I= \langle p_1,\cdots,p_\ell \rangle $ for some $p_1,\cdots,p_\ell\in I$.
\end{thm}

The {\it variety} of an ideal $I$ is the set  of all common {\bf complex}
zeros of the polynomials in $I$:
$$ V(I) \quad = \quad \{X\in \cpx^m: p(X)=0 \,\,\, \hbox{for all} \,\, p \in I\}.$$
The subset of all real points in $V(I)$ is the {\it real variety} of $I$. It is denoted
$$ V^\re(I) \quad = \quad \{X\in \re^m: p(X)=0 \,\,\,
\hbox{for all} \,\, p  \in I\}.$$
If $\,I = \langle p_1, \ldots , p_r \rangle \,$ then
$\,V(I)\,=\,V(p_1,\ldots, p_r)=\{X\in \cpx^m: p_1(X)=\cdots=p_r(X)=0\}$.
An ideal $I \subseteq \re[X]$ is {\it zero-dimensional}
if its variety $V(I)$ is a finite set. This condition
is much stronger than requiring that the real variety
$V^\re(I)$ be a finite set. For example,
$I = \langle X_1^2 + X_2^2 \rangle$ is not zero-dimensional, however
the real variety $V^{\re}(I) = \{(0,0)\}$ consists of one point of the
curve $V(I)$.


A variety $V\subseteq \cpx^m$ is {\it irreducible}
if there do not exist two proper subvarieties $V_1,V_2 \subsetneqq  V$ such that
$V=V_1 \cup V_2$.
The reader should note that in this paper,
``irreducible" means that the set of {\bf complex} zeros cannot be written as
a proper union of subvarieties defined by {\bf real} polynomials.
Given a variety $V \subseteq \cpx^m$, the set of
all polynomials that vanish on $V$ is an ideal
$$\,I(V) \quad = \quad \{p\in \re[x]\,:\, p(u)=0 \,\,\,\hbox{for all} \,\,  u \in V\}. $$
Given any ideal $I$ of $\re[x]$, its {\it radical} is the ideal
$$ \sqrt{I} \quad = \quad  \bigl\{q\in \re[X] \,\,:\,\, q^\ell \in I \,
\mbox{ for some } \ell \in \N\bigr\}.$$
Note that $I \subseteq \sqrt{I}$. We say that $I$ is a {\it radical ideal} if $\sqrt{I}=I$.
Clearly, the ideal $I(V)$ defined by a variety $V$ is a radical
ideal. The following theorems offer a converse to this observation:

\begin{thm} [Hilbert's Weak Nullstellensatz] \label{WeakNull}
\hfill \break
If $I$ is an ideal in $\re[x]$ such that $V(I)=\emptyset$
then $1\in I.$
\end{thm}

\begin{thm} [Hilbert's Strong Nullstellensatz] \label{StrNull}
\hfill \break
If $I$ is an ideal in $\re[x]$ then
$I(V(I)) = \sqrt{I}.$
\end{thm}


\vspace{.3cm}

\noindent
{\it Remark.}  Theorems \ref{WeakNull}, \ref{StrNull} are normally stated
for ideals in $\cpx[x]$.  However, keeping in mind that $V(I)$ lies in
$\cpx^m$, they hold as stated.
\vspace{.2cm}

In real algebraic geometry, we are also interested in
subsets of $\re^m$ of the form
$$ S \quad = \quad \bigl\{\,X\in \re^m \, :\,
p_1(X)=\cdots=p_r(X)=0,
q_1(X)\geq 0, \cdots,q_\ell(X)\geq 0 \bigr\} , $$
where $p_i, q_j \in \re[x]$. Such $S$ is called
a {\it basic closed  semialgebraic set}.
Given $S$ as above, the preorder and linear cones associated with $S$ are
defined as
\begin{align*}
P(S)\,\, &=\,\,\left\{
\left.
\sum_{\theta\in \{0,1\}^\ell}
\sig_{\theta}(X) q_1^{\theta_1}(X)\cdots q_\ell^{\theta_\ell}(X)
\right|
\sig_0,\cdots,\sig_\ell \mbox{ are SOS} \right\} +
\Big\langle p_1,\cdots, p_r \Big\rangle \\
M(S)\,\, &=\,\,\left\{
\left. \sig_0(X)+\sum_{j=1}^\ell q_j(X)\sig_j(X)
\right|
\sig_0,\sig_1,\cdots,\sig_\ell \mbox{ are SOS}
 \right\}+
\Big\langle p_1,\cdots, p_r \Big\rangle .
\end{align*}


A linear cone or preorder $M$ is {\it archimedean} if there exists
$\rho(x)\in M$ such that
the set $\{X\in \re^m: \rho(X) \geq 0\}$ is compact, equivalently, if
there exists $N \in \N$ such that $N - \sum_{i=1}^m x_i^2 \in M$.  Note
that if $M(S)$ or $P(S)$ is archimedean, then  $S$ is compact.

\begin{thm}[Putinar, \cite{putinar}] \label{putinar}
Suppose $M(S)$ is archimedean,
then every polynomial $p(x)$ which is positive
on $S$ belongs to $M(S)$.
\end{thm}

\vspace{.3cm}
\noindent
{\it Remark.} \label{archimedean}
There are examples of compact $S$ for which $M(S)$ is not archimedean
and the conclusion of Putinar's Theorem does not hold.
In the case of the preorder $P(S)$, it is a deep theorem of Schm\"udgen
\cite{Sch} that if $S$ is compact then $P(S)$ is archimedean and any
polynomial which is positive on $S$ is in $P(S)$.  For this reason,
the SOS relaxations $p^*_N$ always converge to the minimum if
$S$ is compact, however,
the relaxations $f^*_N$ may not converge to the minimum.
On the other hand, it is sometimes the case in practice that we know or can
compute some $N \in \N$ such that our semialgebraic set $S$ is contained
in the sphere $\{ N - \sum_{i=1}^m x_i^2 \geq  0 \}$.  In this case, we can
simply add one additional constraint, namely  $N - \sum_{i=1}^m x_i^2 \geq 0$,
and force $M(S)$ to be archimedean.
\vspace{.2cm}

The sets $P(S)$ and $M(S)$ contain the ideal
 $J=\langle p_1,\cdots,p_r\rangle$.
If $J$ is radical and $V(J)$ is finite, we have the following theorem:
\begin{thm}[Parrilo, \cite{Par02}] \label{par_finrad}
Let $S$ and $J$ be defined as in the above.
Suppose $J$ is a zero-dimensional radical
ideal in $\re[X]$. Then a polynomial
$w(x)\in \re[X]$ is nonnegative on S
if and only if $w(x)\in M(S)$.
\end{thm}

For a semialgebraic set, there is a well-known generalization
of the {\it Hilbert's Weak Nullstellensatz}, see e.g. \cite[4.2.13]{PD}.


\begin{thm} \label{posnull}
Suppose $S$ and $P(S)$ are defined as above, then
$S=\emptyset$ if and only if $-1\in P(S)$.
\end{thm}

We need the following lemma from \cite{njw_grad}:

\begin{lemma}[Lemma 3.2,\cite{njw_grad}] \label{lagpoly}
Let $V_1,\cdots,V_r$ be pairwise disjoint varieties of $\cpx^m$.
Then there exist polynomials $p_1,\cdots,p_r\in \cpx[X]$ such that
$p_i(V_j)=\delta_{ij}$, where $\delta_{ij}$ is the Kronecker
delta function.
\end{lemma}
\noindent
Furthermore, if each $V_{\ell}$ is conjugate symmetric,
i.e., a point $z \in \cpx^m$ belongs to $V_{\ell}$ if and only if
its complex conjugate $\bar z \in V_{\ell}$,
then the polynomials $p_{\ell}$ can be chosen
such that $p_{\ell}\in \re[X]$,
since we can replace $p_i(x)$ by $(p_i(X)+\bar p_i(X))/2$,
where $\bar p_i(X)$ is obtained from $p_i(X)$ by
conjugating its coefficients.

\Section{Representations in $P_{KKT}$ and $M_{KKT}$}


In \cite{njw_grad}, it is
shown that if a polynomial $f(x) \in \re[x]$ is globally nonnegative
and its gradient ideal is radical, then $f(x)$ has a representation
as a sum of squares modulo the gradient ideal.
In this section we generalize this result to real polynomials which
are nonnegative on the semialgebraic set $V_{KKT}$:  We will show that
such polynomials have a representation in $P_{KKT}$ modulo the ideal
$I_{KKT}$,  if the later is radical.  Furthermore, in some cases we
can replace the preorder cone $P_{KKT}$ by the linear cone $M_{KKT}$.

Throughout this section we fix a polynomial $f(x) \in \re[x]$ along
with an optimization of the form (1.1)-(1.3) and the corresponding
ideal $I_{KKT}$, variety $V_{KKT}$, the preorder cone $P_{KKT}$ and
the linear cone $M_{KKT}$.


>From Theorem~\ref{par_finrad}, we immediately obtain the following
representation theorem:
\begin{thm}
Assume $I_{KKT}$ is zero-dimensional and radical.
If $f(x)$ is nonnegative on $V_{KKT}^\re \cap {\cal H}$,
then $f(x)$ belongs to $M_{KKT}$.
\end{thm}

Using a proof similar to that of Theorem~3.1 in \cite{njw_grad},
we can remove the restrictive hypothesis that $I_{KKT}$ be zero-dimensional,
however to obtain the most general result we must replace the
linear cone $M_{KKT}$ by the preorder cone $P_{KKT}$.

\begin{thm} \label{radrepr}
Assume $I_{KKT}$ is radical.
If $f(x)$ is nonnegative on $V_{KKT}^\re \cap {\cal H}$,
then $f(x)$ belongs to $P_{KKT}$.
\end{thm}

We need a generalization of a lemma from \cite{njw_grad}:

\begin{lemma} \label{conval}
Let $W$ be an irreducible component of $V_{KKT}$.
Then $f(x)$ is constant on $W$.
\end{lemma}
%

\proof
We first note that
$$
F(x) = f(x) + \sum_{i=1}^s \lmd_i g_i(x) + \sum_{j=1}^t \nu_j h_j(x)
$$
is equal to $f(x)$ on $V_{KKT}$,
and the right hand side has zero gradient on $V_{KKT}$. With this in
mind, the proof of \cite[3.3]{njw_grad} generalizes easily to
this case.
\eproof

\medskip

\noindent {\it Proof of Theorem \ref{radrepr}.} \quad
Decompose $V_{KKT}$ into its irreducible components, then
by Lemma~\ref{conval}, $f(x)$ is constant on each of them.
Let $W_0$ be the union of all the components whose intersection
with ${\cal H}$ is empty, and group together the components
on which $f(x)$ attains the same value, say  $W_1, \dots , W_r$.
Suppose $f(x) = \alpha_i \geq 0$ on $W_i$.

We have $V_{KKT}=W_0 \cup W_1 \cup \cdots \cup W_r$,
and $W_i$ are pairwise disjoint. Note that by our definition of irreducible,
each $W_i$ is conjugate symmetric.
By Lemma~\ref{lagpoly}, there exist polynomials
$p_0,p_1,\cdots,p_r\in \re[x,\lmd,\nu]$ such that
$p_i(W_j)=\delta_{ij}$, where $\delta_{ij}$
is the Kronecker delta function.

By assumption,  $W_0 \cap {\cal H}=\emptyset$ and so, by
Theorem~\ref{posnull}, there are
SOS polynomials $v_{\theta}\,\ (\theta \in \{0,1\}^t)$
such that
$$
-1 \equiv   \sum_{\theta \in \{0,1\}^t}
v_{\theta} h_1^{\theta_1}\cdots h_t^{\theta_t}
\overset{def}{=}v_0
\,\,\,\  mod \,\,\,\ I(W_0).
$$
We have $f=(f +\frac{1}{2})^2-(f^2+(\frac{1}{2})^2)=f_1+v_0\cdot f_2$
for the SOS polynomials $f_1 = (f + \frac{1}{2})^2, f_2 = f^2 + (\frac12)^2$.
Then
$$f \equiv f_1+v_0f_2 \equiv
\sum_{\theta \in \{0,1\}^t}
u_{\theta} h_1^{\theta_1}\cdots h_t^{\theta_t}
\overset{def}{=} q_0
\,\,\,\  mod \,\,\,\ I(W_0)
$$
for some SOS polynomials $u_{\theta}\,\ (\theta \in \{0,1\}^t) $.
Recall that $f(x) = \alpha_i$, a constant, on each $W_i(1\leq i \leq r)$.
Set $q_i(x) = \sqrt{\alpha_i}$, then $f(x)=q_i(x)^2$ on $I(W_i)$.

Now let $q=q_0 (p_0)^2+ \sum_{i=1}^r (q_i p_i)^2$.
Then $f - q$ vanishes on $V_{KKT}$ and hence $f - q \in I_{KKT}$ since
$I_{KKT}$ is radical.  It follows that $f \in P_{KKT}$.
\eproof

\vspace{.3cm}
\noindent
{\it Remark.} The assumption that
$I_{KKT}$ is radical
is needed in
Theorem~\ref{radrepr}, as shown by
Example~3.4 in \cite{njw_grad}.
However, when $I_{KKT}$ is not radical, the conclusion also holds
if $f(x)$ is strictly positive on $V_{KKT}^\re$.
\vspace{.2cm}

\begin{thm} \label{poskkt}
If $f(x)$ is strictly positive on $V_{KKT}^\re \cap {\cal H}$
then
$f(x)$ belongs to $P_{KKT}$.
\end{thm}
\proof
As in the proof of Theorem \ref{radrepr}, we decompose $V_{KKT}$ into
subvarieties $W_0,W_1,\cdots,W_r$ such that $W_0 \cap {\cal H} =
\emptyset$, and for $i = 1, \dots r$, $W_i \cap {\cal H} \neq \emptyset$ and
$f$ is constant on $W_i$.
Since each $W_i$, $i>0$ contains at least one real point and
$f(x)>0$ on $V_{KKT}^\re$, each $\alpha_i>0$.  The $W_i$ were chosen
so that each $\alpha_i$ is distinct, hence
the $W_i$'s are pairwise disjoint.

Consider the primary decomposition  $I_{KKT}= \cap_{i=0}^r J_i $
corresponding to our decomposition of $V_{KKT}$, i.e.,
$V(J_i)=W_i$ for $i=0,1,\cdots,r$.
Since $W_i \cap W_j =\emptyset$, we have
$J_i+J_j=\re[x,\lmd,\nu]$ by Theorem~\ref{WeakNull}.
The Chinese Remainder Theorem, see e.g. \cite[2.13]{Eisenbud}, implies that
there is an isomorphism
$$
\rho: \re[x,\lmd,\nu]\big/ I_{KKT} \rightarrow
\re[x,\lmd,\nu]\big/ J_0 \times
\re[x,\lmd,\nu]\big/ J_1 \times \cdots \times
\re[x,\lmd,\nu]\big/ J_r.
$$

For any $p\in \re[x,\lmd,\nu]$, let
$[p]$ and $\rho([p])_i$ denote the equivalence classes of
$p$ in $\re[x,\lmd,\nu]\big/ I_{KKT}$
and $\re[x,\lmd,\nu]\big/ J_i$ respectively.

Recall that that $V(J_0) \cap {\cal H}=\emptyset$, hence by
Theorem~\ref{posnull}
there exist SOS polynomials
$u_{\theta}\,\ (\theta \in \{0,1\}^t)$ such that
$$ -1 \equiv  \sum_{\theta \in \{0,1\}^t}
u_{\theta} \rho([h_1^{\theta_1}])_0 \cdots \rho([h_t^{\theta_t}])_0
\overset{def}{=}  u_0
\,\,\,\ \mbox{ mod }\,\,\,\ J_0\ .
$$
As in the proof of Theorem \ref{radrepr}, we write
$f =    f_1  -  f_2$
for SOS polynomials $  f_1, f_2 $ and then we have
$$
f \equiv f_1  + u_0  f_2
\equiv
\sum_{\theta \in \{0,1\}^t}
v_{\theta} (\rho([h_1^{\theta_1}]))_0 \cdots (\rho([h_t^{\theta_t}]))_0
\overset{def}{=} q_0 \thinspace \mbox{ mod } J_0
$$
for some SOS polynomials $v_{\theta}\,\ (\theta \in \{0,1\}^t) $.
Thus the preimage $\rho^{-1}((q_0,0,\cdots,0)) \in P_{KKT}$.

Now on each $W_i$, $1\leq i\leq r$, $f(x) = \alpha_i>0$,
and hence $(f(x)\big/\alpha_i) - 1$ vanishes on $W_i$.
Then by Theorem \ref{StrNull} there is
$\ell \in \N$ such that
$(f(x)\big/\alpha_i-1)^\ell \in J_i$.
>From the binomial theorem, it follows that
$$
(1+(f(x)\big/\alpha_i-1))^{1/2} \equiv \sum_{k=1}^{\ell-1}
\binom{1/2}{k} (f(x)\big/\alpha_i-1)^k \overset{def}{=} q_i\big/ \sqrt{\alpha_i}
\,\,\,\ \mbox{ mod } J_i \,\ .
$$
Thus $(\rho([f]))_i= q_i^2 $ is SOS in $\re[x,\lmd,\nu]\big/ J_i$, and hence
$\rho^{-1}(q_i^2 e_{i+1})$ is SOS in $\re[x,\lmd,\nu]\big/ I_{KKT}$,
where $e_{i+1}$ is the $(i+1)$-st standard unit vector in $\re^{r+1}$.

Finally, we see that $\rho([f]) =
(q_0,q_1^2,\cdots,q_r^2)$. The preimage of the latter is
$$
\rho^{-1}\big((q_0,q_1^2,\cdots,q_r^2)\big) =
\rho^{-1}\big(q_0 e_1)\big) + \sum_{i=1}^r\rho^{-1}\big(q_i^2 e_{i+1}\big),
$$
which implies that $f\in P_{KKT}$.
\eproof

\vspace{.3cm}

\noindent
{\it Remark.} The conclusions in Theorem~\ref{radrepr}
and Theorem~\ref{poskkt} can not be strengthened to show that
$f(x) \in M_{KKT}$. The following is a counterexample.
\vspace{.2cm}

\begin{exm} \label{podnelin}
Consider the optimization
\begin{align*}
\min & \,\,\,\  f(x)=(x_3-x_1^2x_2)^2-1+\eps \\
s.t. & \,\,\,\  h_1(x)=1-x_1^2\geq 0 \\
& \,\,\,\ h_2(x)=x_2\geq 0 \\
& \,\,\,\ h_3(x)=x_3-x_2-1\geq 0
\end{align*}
where $0< \eps < 1$.
>From the constraints, we can easily observe that
the global minimum $f^*=\eps>0$ which is
attained at $x^*=(0,0,1)$.
Its KKT ideal
\begin{align*}
I_{KKT}=\Big\langle
&2x_1x_2(x_3-x_1^2x_2)-\nu_1x_1,
2x_1^2(x_3-x_1^2x_2)+\nu_2-\nu_3, \\
&2(x_3-x_1^2x_2)-\nu_3,
 \nu_1(1-x_1^2),\nu_2x_2,\nu_3(x_3-x_2-1)  \Big\rangle
\end{align*}
is radical (verified in Macaulay 2 \cite{Macaulay}).
However, we can not find
SOS polynomials $\sig_0,\sig_1,\sig_2,\sig_3$ and general polynomials
$\phi_1,\phi_2,\phi_3$ such that
$$ f(x) = \sig_0 + \sig_1 h_1 +
\sig_2 h_2 + \sig_3 h_3+
\phi_1 (\frac{\pt f}{\pt x_1}-\nu_1x_2) +
\phi_2 (\frac{\pt f}{\pt x_2}-\nu_2+\nu_3)+
\phi_3 (\frac{\pt f}{\pt x_3}-\nu_3).$$
Suppose to the contrary that they exist.
Plugging $\nu=(0,0)$ in the above identity yields
$$0 = 1-\eps+\sig_0 + \sig_1 (1-x_1^2) +
 \sig_2 x_2 +\sig_3 (x_3-x_2-1)+
\phi (x_3-x_1^2x_2)$$
where
$\phi=-4x_1\phi_1-x_1^2\phi_2+2\phi_3-(x_3-x_1^2x_2)$.
Now substitute $x_3=x_1^2x_2$ in the above, yielding
$$\sig_3 ((1-x_1^2)x_2+1) = 1-\eps+\sig_0 + \sig_1 (1-x_1^2) +
 \sig_2 x_2.$$
Here $\sig_0,\sig_1,\sig_2,\sig_3$ are now considered as
SOS polynomials in $(x_1,x_2)$.
Since $1-\eps>0$, $\sig_3$ can not be the zero polynomial.
If $\sig_3=\sig_3(x_1)$ is independent of $x_2$,
we can derive a contradiction using an argument identical
to the argument in the proof of
of \cite[Thm. 2]{Powers_rectangle}.
Thus $2m=\mbox{deg}_{x_2}\sig_3(x_1,x_2) \geq 2$
and $2d=\mbox{deg}_{x_1}\sig_3(x_1,x_2) \geq 0$.
On the left hand side, the leading term is of the form
$A\cdot x_1^{2d+2}x_2^{2m+1}$ with coefficient
$A<0$. Since the degree in $x_2$ on the left hand side is odd,
the leading term on the right hand side must come
from $\sig_2(x_1,x_2)x_2$,
and is of the form like $B\cdot x_1^{2d}x_2^{2m+1}$
with $B>0$. This is a contradiction.
Therefore we can conclude that $f(x) \notin M_{KKT}$.
%
\end{exm}

\Section{ Convergence of the Lower Bounds}

In this section, we will show that the lower bounds
$\{p^*_N\}$ obtained from \reff{sospod1}-\reff{sospod2}
converge to $f^*$ in \reff{opt1}-\reff{opt3}.
The conclusions in Section~4 of \cite{njw_grad}
can generalized, based on Theorem~\ref{radrepr}
and Theorem~\ref{poskkt} in the preceding section.
However, we need an extra assumption to ensure
the convergence of $\{f^*_N\}$.

\begin{thm} \label{convergence}
Assume $f^*$ is finite and
the global optimizers $x^*$ of \reff{opt1}-\reff{opt3}
satisfy the KKT system \reff{kkt1}-\reff{kkt3}.
Then $\underset{N\to \infty}{\lim}p^*_N=f^*$.
Furthermore, if $I_{KKT}$ is radical,
then there exists some $N\in \N$ such that
$p^*_N=f^*$, i.e.,
the SOS relaxations \reff{sospod1}-\reff{sospod2}
converge in finitely many steps.
\end{thm}
\proof
The sequence $\{p^*_N\}$ is monotonically increasing,
and $p^*_N \leq f^*$ for all $N\in \N$, since $f^*$ is
attained by $f(x)$ in the KKT system \reff{kkt1}-\reff{kkt3}
by assumption and the constraint
\reff{sosrlx2} implies that $\gamma \leq f^*$.
Now for arbitrary $\eps >0$, let $\gamma_\eps=f^*-\eps$
and replace $f(x)$ by $f(x)-\gamma_\eps$ in \reff{opt1}-\reff{opt3}.
The KKT system remains unchanged, and $f(x)-\gamma_\eps$ is strictly
positive on $V_{KKT}^\re$.
By Theorem~\ref{poskkt}, $f(x)-\gamma_\eps \in P_{KKT}$.
Since $f(x)-\gamma_\eps$ is fixed, there must exist some integer $N_1$
such that $f(x)-\gamma_\eps \in P_{N_1,KKT}$. Hence
$ f^*-\eps \leq p^*_{N_1} \leq f^*$.
Therefore we have that $\underset{N\to \infty}{\lim}p^*_N=f^*$.

Now assume that $I_{KKT}$ is radical.
Replace $f(x)$ by $f(x)-f^*$ in \reff{opt1}-\reff{opt3}.
The KKT system still remains the same, and $f(x)-f^*$ is
now nonnegative on $V_{KKT}^\re$.
By Theorem~\ref{radrepr}, $f(x)-f^* \in P_{KKT}$.
So there exists some integer $N_2$
such that $f(x)-f^* \in P_{N_2,KKT}$, and hence
$P^*_{N_2} \geq f^*$. Then $p^*_{N}\leq f^*$ for all $N$
implies that $p^*_{N_2}=f^*$.
\eproof

\vspace{.3cm}

\noindent
{\it Remarks.}
(1) In Lasserre's method \cite{lasserre}, a
sequence of lower bounds that converge to $f^*$ asymptotically
can be obtained when the feasible region ${\cal F}$ is compact;
but those lower bounds usually do not converge in finitely many steps.
However, from Theorem~\ref{convergence}, we see that when $I_{KKT}$
is radical then
the lower bounds $\{p^*_N\}$ converge in finitely many steps, even
if ${\cal F}$ is not compact.
This implies that the lower bounds $\{p^*_N\}$ may have better
convergence even in the case where ${\cal F}$ is compact.

(2) The assumption in Theorem~\ref{convergence}
can not be removed, which is illustrated by the
following example.

\vspace{.2cm}

\begin{exm} \label{kktfails}
Consider the optimization:
$ \min\,\ x \mbox{ s.t. } x^3 \geq 0. $
Obviously $f^*=0$ and the global minimizer $x^*=0$.
However, the KKT system
$$
1 - \nu \cdot 3x^2=0,\,\,\
\nu \cdot x^3=0,\,\,\ x^3\geq0,\,\,\ \nu\geq 0
$$
is not satisfied, since $V_{KKT}=\emptyset$.
Actually we can see that the lower bounds $\{f^*_N\}$ given by
\reff{sosrlx1}-\reff{sosrlx2} tend to infinity.
By Theorem~\ref{WeakNull},
$V_{KKT}=\emptyset$ implies that $1\in P_{KKT}$, i.e.,
$$
(1+3\nu x^2)(1-3\nu x^2)+9\nu^2 x \cdot \nu x^3=1.
$$
In the SOS relaxation \reff{sosrlx1}-\reff{sosrlx2}, for arbitrarily large
$\gamma$, $x-\gamma \in P_{KKT}$, since
$$
x-\gamma =(x-\gamma)(1+3\nu x^2)(1-3\nu x^2)+9\nu^2 x (x-\gamma) \cdot \nu x^3
\in P_{KKT}.
$$
Thus $p^*_8=\infty$.
In this example, the conclusion in Theorem~\ref{convergence} does not hold.
\end{exm}

The convergence of lowers bounds
$\{f^*_N\}$ cannot be guaranteed,
as we see in
Example~\ref{podnelin}.
In that example,
replace the objective by the perfect square
$(x_3-x_1^2x_2)^2$.
Then $f^*=1$, but we do not have
$\underset{N\to \infty}{\lim}f^*_N=1$.
>From the arguments there,
we can see that
$f(x)-(1-\eps) \notin M_{KKT}$ for all
$0<\eps <1$, which implies that
$f^*_N \leq 0$.
But $f^*_N \geq 0$ is obvious since
$(x_3-x_1^2x_2)^2$ is a perfect square.
Therefore $\underset{N\to \infty}{\lim}f^*_N=0<1=f^*$, i.e.,
the lower bounds $\{f^*_N\}$ obtained from
\reff{sosrlx1}-\reff{sosrlx2}
may not converge.

On the other hand, the situation is often not that bad in practice.
In the examples in the rest of this paper,
it always happens that
$\underset{N\to \infty}{\lim}p^*_N=\underset{N\to \infty}{\lim}f^*_N=f^*$.
If we further assume that
$M_{KKT}$ is archimedean
then it must hold that
$\underset{N\to \infty}{\lim}p^*_N=\underset{N\to \infty}{\lim}f^*_N=f^*$
from Theorem~\ref{putinar} (Putinar).
This is the generalization of assumption~4.1 in \cite{lasserre}.
See also the remark after Theorem~\ref{putinar}.

\vspace{.2cm}

The SOS relaxation \reff{sosrlx1}-\reff{sosrlx2}
can be solved using software SOSTOOLS \cite{soostools}.
The dual problem of \reff{sosrlx1}-\reff{sosrlx2}
is to minimize a linear functional
over some linear moment matrix inequalities.
It can also be obtained by applying
moment matrix methods to minimize $f(x)$
over the semialgebraic set defined by KKT system
\reff{kkt1}-\reff{kkt3}.
The dual problem can be solved using
software Gloptipoly \cite{Gloptipoly}.
Actually, the formulations of SOS relaxations and
moment matrix methods are
dual to each other, see \cite{lasserre, Laurent}.
The SOS relaxations
\reff{sosrlx1}-\reff{sosrlx2}
not only give the lower bounds $f^*_N$,
but also the information about
global minimizers $x^*$ and their
Lagrange multipliers $(\lmd^*,\nu^*)$.
SOSTOOLS can extract the minimizer
if the moment matrix has rank one.
Gloptipoly can
also find the lower bounds, and
extract (\cite{HenLas}) the global minimizers
when the moment matrix satisfies some rank condition.
Gloptipoly does not need the moment matrix to be rank one.
The tricks to extract global minimizers in
Section~5.2 in \cite{njw_grad}
can be applied here directly to find $(x^*,\lmd^*,\nu^*)$,
so omit further discussion.
For more details about how extracting
minimizers from SOS relaxations or
moment matrix methods, see \cite{HenLas}.

\begin{exm}[Exercise 2.18, \cite{pardalos_intr_globopt}]
Consider the global optimization:
\begin{align*}
\min &\,\,\,\ (-4x_1^2+x_2^2)(3x_1+4x_2-12) \\
s.t. &\,\,\,\  3x_1-4x_2\leq 12,\,\,\ 2x_1-x_2\leq 0,\,\,\ -2x_1-x_2 \geq 0.
\end{align*}
The global minimum $f^*=-18.6182$ and the minimizer
$x^*=(-24/55,128/55)\approx (-0.4364, 2.3273)$.
The lower bound obtained from \reff{sosrlx1}-\reff{sosrlx2}
is $f^*_4=-18.6182$.
The extracted minimizer $\hat x=(-0.4364, 2.3273)$.
\end{exm}


\begin{exm} Consider the Quadratically Constrained
Quadratic Program (QCQP):
\begin{align*}
\min &\,\,\,\ -\frac{4}{3}x_1^2+\frac{2}{3}x_2^2-2x_1x_2 \\
s.t. &\,\,\,\  x_2^2-x_1^2\geq 0,\,\,\ -x_1x_2\geq 0.
\end{align*}
The global minimum $f^*=0$ and minimizer $x^*=(0,0)$.
The feasible region ${\cal F}$ defined by the constraints is
non-compact.
The lower bound returned by
\reff{sosrlx1}-\reff{sosrlx2}
is $f^*_4=-2.6 \times 10^{-15}$
(Note: this computation was done in double precision
floating point, with round off error bounded by
$2^{-53}\sim 10^{-16}$).
The extracted minimizer is
$\hat x=(6.1\times 10^{-16}, -9.0\times 10^{-17})$
and the Lagrange multiplier is
$\hat \nu=(0.3884, 0.3909)$.
\end{exm}

\Section{Optimization over Some Special Semialgebraic Sets}

In problem \reff{sosrlx1}-\reff{sosrlx2}, the polynomials
are in $(x,\lmd,\nu)\in \re^{n+s+t}$ which means that  when there are
many constraints, the problem is very expensive to solve.
If $u(x,\lmd,\nu)$ is a polynomial of degree $d$,
it can have $\binom{n+s+t+d}{d}$ coefficients;
this will be huge for large $s$, $t$, or $d$.
Frequently, if the polynomials $g_i(x)$ and $h_j(x)$ are
of some special form, then the
KKT system \reff{kkt1}-\reff{kkt3} can be simplified
and hence the SOS relaxations
\reff{sosrlx1}-\reff{sosrlx2} will be easier to solve.
In this section we look at the case where
$\{x\in \re^n: h_1(x),\cdots,h_t(x)\geq 0\}$
is the nonnegative orthant $\re^n_{+}$ or
the box $[a,b]_n$ and show how these type of problems can be
simplified.

\subsection{Minimizing Over the Nonnegative Orthant $\re^n_+$}

In this subsection, suppose the inequality constraints \reff{opt3} are
the standard constraints for the
nonnegative orthant
$\re^n_+:=\{x\in \re^n: x_1\geq 0,\cdots, x_n\geq 0\}$.
The constraints are
of the form
\begin{align*}
g_1(x)=\cdots=g_s(x)=0,\,\,\,\  x \in \re_+^n.
\end{align*}
Then the KKT system \reff{kkt1}-\reff{kkt3} becomes
\begin{align*}
\nabla f(x) &+\sum_{i=1}^s \lmd_i \nabla g_i(x) - \nu =0, \\
g_1(x)&=\cdots=g_s(x)=0,  \\
 x_k\nu_k &=0,\,\ k=1,\cdots,n, \\
x&\in \re^n_+,\,\,\,\  \nu \in \re^n.
\end{align*}
In this KKT system, the variable $\nu$ can be solved for explicitly.
By eliminating $\nu$, the above system simplifies to
\begin{align}
x_k(\frac{\pt f}{\pt x_k} &+
\sum_{i=1}^s\lmd_i\frac{\pt g_i}{\pt x_k} )
=0,\,\ k=1,\cdots,n \label{kktorthant1}\\
g_1(x)=&\cdots=g_s(x)=0. \label{kktorthant2}
\end{align}
We define cones $M_{KKT}^{\re^n_+}$ and $M_{N,KKT}^{\re^n_+}$
similar to the definition of
$M_{KKT}$ and $M_{N,KKT}$ (see Section~1),
define
associated to the above simplified system.
Note that $M_{KKT}^{\re^n_+},M_{KKT}^{\re^n_+}\subset \re[x,\lmd]$ and
the Lagrange multiplier $\nu$ does not appear.
Similar to \reff{sosrlx1}-\reff{sosrlx2},
a sequence $\{\hat f^*_N\}$ of lower bounds of \reff{opt1}-\reff{opt3}
can be obtained by the following SOS relaxations:
\begin{align}
\hat f^*_N= \max_{\gamma \in \re} &\,\,\,\ \gamma  \label{sosorthant1} \\
s.t. &\,\,\,\ f(x)-\gamma \in M_{N,KKT}^{\re^n_+}. \label{sosorthant2}
\end{align}
Now the indeterminates in the above SOS program
are $(x,\lmd)$ instead of $(x,\lmd,\nu)$.
Thus a polynomial $u(x,\lmd)$ of degree $d$
has at most $\binom{n+s+d}{d}$ coefficients,
which is much smaller than $\binom{n+s+t+d}{d}$
when $t$ is large.
This makes solving \reff{sosorthant1}-\reff{sosorthant2}
much less expensive.

Since $\nu=(\nu_1,\cdots,\nu_t)$ are eliminated
by direct substitutions,
systems \reff{kkt1}-\reff{kkt3} and
\reff{kktorthant1}-\reff{kktorthant2} are
equivalent. Thus we see that
$f(x)-\gamma \in M_{N_1,KKT}$
if and only if
$f(x)-\gamma \in M_{N_2,KKT}^{\re^n_+}$,
for some integers $N_1$ and $N_2$.
Therefore the lower bounds $\{\hat f^*_N\}$
have the same property of convergence
as $\{f^*_N\}$ obtained
from \reff{sosrlx1}-\reff{sosrlx2}.

If, in addition, the equality constraints \reff{opt2}
are hyperplanes, i.e., the constraints are
the standard simplex:
\begin{align*}
Ax=b, \,\,\,\
x \geq 0
\end{align*}
where $A\in \re^{s \times n}, b\in \re^s$,
then the KKT system \reff{kkt1}-\reff{kkt3}
can be reduced to
\begin{align*}
x_k(\frac{\pt f}{\pt x_k} & +a_k^T\lmd) =0,\,\ k=1,\cdots,n \\
Ax&=b,\,\,\, x\geq 0
\end{align*}
where $a_k\in \re^s$ is the $k$-th column of matrix $A$.

Furthermore, if $Ax=b$ consists of a single equation
$a^Tx=b\ne 0$,
then $\lmd=-\frac{x^T\nabla f(x)}{b}$ and
the KKT system has the simpler form
\begin{align*}
x_k(\frac{\pt f}{\pt x_k} & -\alpha_k \frac{x^T\nabla f(x)}{b} ) =0,\,\ k=1,\cdots,n \\
a^Tx&=b,\,\,\, x\geq 0
\end{align*}
where $a=[\alpha_1,\cdots,\alpha_n]^T$.

Based on the above two simplified KKT systems,
SOS relaxations similar to \reff{sosorthant1}-\reff{sosorthant2}
can be obtained immediately, improving the computational efficiency.

\begin{exm}
[Test Problem~2.9, \cite{pardalos_test_prob}]
Consider the Maximum Clique Problem for $n=5$:
\begin{align*}
\min &\,\,\,\ -\big(\sum_{i=1}^4 x_ix_{i+1} + x_1x_5+x_1x_4+x_2x_5+x_3x_5\big) \\
s.t. &\,\,\,\ x_1+x_2+x_3+x_4+x_5 =1 \\
&\,\,\,\ x_1,x_2,x_3,x_4,x_5 \geq 0.
\end{align*}
The global minimum $f^*=-1/3$ and minimizers $x^*$ are
$(1/3,1/3,0,0,1/3)$,
$(1/3,0,0,1/3,1/3)$,
$(0,1/3,1/3,0,1/3)$, and
$(0,0,1/3,1/3,1/3)$.
The lower bound obtained from
\reff{sosorthant1}-\reff{sosorthant2}  is
$\hat f^*_4=-0.33333333378814$.
The difference $f^*-\hat f^*_4\approx 4.5\times 10^{-10}$.
\end{exm}

\begin{exm}[Exercise~1.20, \cite{pardalos_intr_globopt}]
Consider the optimization:
\begin{align*}
\min &\,\,\,\ \sum_{i=1}^{n-1} x_i^2x_{i+1}+x_n^2x_1 \\
s.t. &\,\,\,\ \sum_{x_i=1}^n x_i =1,\,\,\,\ x\geq 0.
\end{align*}
The global minimum $f^*=0$ and
the minimizers are the vertices of the simplex
defined by the constraints.
The lower bound obtained from \reff{sosorthant1}-\reff{sosorthant2}
is $\hat f_4^*=-4.0 \cdot 10^{-8}$.
\end{exm}

\begin{exm}
$f(x)=x^THx$ and the constraints are
$0\leq x \leq e$, where
$x\in \re^5$ and $e=[1,1,1,1,1]^T$,
and
$$
H=\bbm
1 & -1 &  1  & 1  & -1\\
 -1  & 1  & -1  & 1  & 1\\
 1  & -1 &  1 &  -1 &  1 \\
1  & 1  & -1  & 1  & -1 \\
-1  & 1  & 1 &  -1  & 1
\ebm
$$
is a co-positive matrix (\cite{Par01, Pths}),
i.e., $f(x)\geq 0\,\,\,\forall x\geq 0$.
If each $x_i$ is replaced by $x_i^2$,
then the resulting quartic polynomial
is nonnegative, but not SOS.
Consider the Quadratic Program (QP):
\begin{align*}
\min &\,\,\,\,\ x^THx \\
s.t. &\,\,\,\,\  x_1,x_2,x_3, x_4,x_5 \geq 0.
\end{align*}
The lower bound obtained from \reff{sosorthant1}-\reff{sosorthant2} is
$\hat f^*_2=-3.35\times 10^{-9}$.
Actually, we have the following decomposition
$$
x^THx= 0 + \sum_{i=1}^5 2\cdot ( x_i \cdot h_i^T x)
$$
in \reff{sosorthant1}-\reff{sosorthant2}.
Here $h_i$ is the $i$-th column of matrix $H$.
\end{exm}

\subsection{Minimizing Over the Box}

In this subsection,
we consider the case that \reff{opt3} are given by box constraints,
i.e., $ x\in [a,\,\ b]_n$ where
$a=\bbm a_1 & \cdots & a_n\ebm^T$ and
$b=\bbm b_1 & \cdots & b_n\ebm^T$.
Here we assume that $a<b$.
In this case, the feasible region ${\cal F}$ is compact,
and Lasserre's method \cite{lasserre}
can be applied here.
However, as remarked after Theorem~\ref{convergence}, if $I_{KKT}$ is
radical then our method will converge after finitely many steps.
Usually Lasserre's method has only asymptotic convergence.

Now the KKT system \reff{kkt1}-\reff{kkt3} has the form
\begin{align*}
\nabla f(x) &+\sum_{i=1}^s \lmd_i \nabla g_i(x) - \nu + \mu =0, \\
g_1(x)&=\cdots=g_s(x)=0, \\
(x_k -a_k)\nu_k &=0,\,\ (b_k-x_k)\mu_k=0,\,\,\ k=1,\cdots,n, \\
x-a &\geq 0,\,\,\ b-x \geq 0,\,\,
\end{align*}
where $\nu_i(\mu_i,\lmd_i)$ is the $i$-th component of
Lagrange multipliers $\nu(\mu,\lmd)$ respectively.
One good property of this KKT system is that
the vectors $\nu$ and $\mu$ can be solved for explicitly.
Eliminating $\nu$ and $\mu$, we obtain
\begin{align*}
(\frac{\pt f}{\pt x_k}&+\sum_{i=1}^s \lmd_i \frac{\pt g_i}{\pt x_k})
(x_k -a_k)(b_k-x_k)=0,\,\,\ k=1,\cdots,n, \\
g_1(x)&=\cdots=g_s(x)=0,\,\ x-a \geq 0,\,\,\ b-x \geq 0.
\end{align*}
Like the definition of
$M_{KKT}^{\re^n_+}$ and $M_{N,KKT}^{\re^n_+}$ (see the preceding subsection),
define the cones $M_{KKT}^{[a,b]_n}$ and $M_{N,KKT}^{[a,b]_n}$
associated with the above simplified KKT system,
where $M_{KKT}^{[a,b]_n},M_{d,KKT}^{[a,b]_n}\subset \re[x,\lmd]$.
Similar to \reff{sosorthant1}-\reff{sosorthant2},
a sequence of lower bounds $\{\tilde f^*_N\}$ of \reff{opt1}-\reff{opt3}
can be obtained by the following SOS relaxations:
\begin{align}
\tilde f^*_N= \max_{\gamma \in \re} &\,\,\,\ \gamma  \label{sosbox1} \\
s.t. &\,\,\,\ f(x)-\gamma \in M_{N,KKT}^{[a,b]_n}. \label{sosbox2}
\end{align}
Now a polynomial $u(x,\lmd)$ of degree $d$ in $M_{N,KKT}^{[a,b]_n}$
has at most $\binom{n+s+d}{d}$ coefficients,
which is much smaller than $\binom{n+s+2n+d}{d}$,
the number of coefficients of one polynomial of degree $d$
in $M_{N,KKT}$.
So \reff{sosbox1}-\reff{sosbox2} can be solved
much more efficiently.
Similarly as $\{\hat f^*_N\}$,  the lower bounds
$\{\tilde f^*_N\}$ have the same properties of convergence
as $\{f^*_N\}$.

Consider the special case that
$f(x)=\frac{1}{2}x^THx+g^Tx$ is a quadratic function
and there are no equality constraints.
Here $g\in \re^n$ and $H=H^T\in \re^{n\times n}$
is symmetric.  The the above KKT system can be further reduced to
\begin{align*}
(h_k^Tx&+g_k)
(x_k -a_k)(b_k-x_k)=0,\,\,\ k=1,\cdots,n, \\
x-&a \geq 0,\,\,\ b-x \geq 0.
\end{align*}
Here $h_k(g_k)$ is the $k$-th row (component) of
arrays $H(g)$.
Finding the global minimum of a general nonconvex quadratic
function over a box is an NP-hard problem.
The relaxations \reff{sosbox1}-\reff{sosbox2} provides a new approach
for such nonconvex quadratic programming.

\begin{exm}
[Test Problem~4.7, \cite{pardalos_test_prob}]
Consider the optimization:
\begin{align*}
\min &\,\,\,\ -12x_1-7x_2+x_2^2\\
s.t. &\,\,\,\ -2x_1^4+2-x_2=0 \\
&\,\,\,\ 0\leq x_1 \leq 2,\,\,\,\ 0\leq x_2 \leq 3.
\end{align*}
The best known objective value is $-16.73889$.
The lower bound obtained from \reff{sosbox1}-\reff{sosbox2}
is $\tilde f^*_6=-16.73889$.
So $f^*=\tilde f^*_6$.
The extracted minimizer
$\tilde x=(0.7175, 1.4698)$
and Lagrange multiplier $\tilde \lmd=-4.0605$.
\end{exm}

\begin{exm}
[Test Problem~2.1, \cite{pardalos_test_prob}]
Consider the optimization:
\begin{align*}
\min &\,\,\,\  42x_1+ 44x_2+ 45x_3+ 47x_4+ 47.5x_5 -50 \sum_{i=1}^5 x_i^2 \\
s.t. &\,\,\,\  20x_1+12x_2+11x_3+7x_4+4x_5 \leq 40 \\
&\,\,\,\  0\leq x_1,x_2,x_3,x_4,x_5 \leq 1  .
\end{align*}
The global minimum $f^*=-17$ and the
minimizer $x^*=(1,1,0,1,0)$.
The lower bound obtained form \reff{sosbox1}-\reff{sosbox2}
is $\tilde f^*_6=-17.00$.
The extracted minimizer
$\tilde x=(1.00,1.00,0.00,1.00,0.00)$ and
Lagrange multiplier $\tilde \nu=0.1799$.

\end{exm}

\begin{exm}
[Exercise~2.22, \cite{pardalos_intr_globopt}]
Consider the Maximum Independent Set Problem
\begin{align*}
\min &\,\,\,\  -\sum_{i=1}^nx_i +\sum_{(i,j)\in E} x_i x_j \\
s.t. &\,\,\,\ 0 \leq x_i \leq 1,\,\,\ i=1,\cdots,n.
\end{align*}
The negative of the global minimum $-f^*$ equals
the cardinality of the maximum independent vertex
set of $G=(V,E)$.
Let $G$ be a pentagon with two diagonals
which do not intersect in the interior.
Now $n=5$ and $f^*=-2$.
The lower bound obtained from \reff{sosbox1}-\reff{sosbox2}
is $\tilde f^*_4=-2.00$.
\end{exm}

\begin{exm}
[Exercise~1.32,  \cite{pardalos_intr_globopt}]
Consider the optimization:
\begin{align*}
\min &\,\,\,\  \prod_{i=1}^nx_i-\sum_{i=1}^nx_i \\
s.t. &\,\,\,\  0\leq a \leq x_1,\cdots,x_n \leq b.
\end{align*}
The global minimum $f^*=a^n-na$ when $a\geq 1$.
For $n=4,a=2,b=3$,
the lower bound obtained from \reff{sosbox1}-\reff{sosbox2}
is $\tilde f^*_6=8.00$. The extracted
minimizer is $\tilde x=(2.00,2.00,2.00,2.00)$.
\end{exm}

\Section{Conclusions}

This paper generalizes most of the theorems
in \cite{njw_grad} from optimizations constrained by algebraic sets
to optimizations constrained
by semialgebraic sets, under the assumption that
the global minimizers satisfy the KKT system.
The special structures of the KKT system
are exploited to accelerate the algorithm
when the constraints include
the nonnegative orthant $\re^n_+$ or
the standard box $[a, b]_n$.

In general, the SOS relaxations \reff{sosrlx1}-\reff{sosrlx2}
are very hard to solve when there are many constraints,
which introduces many Lagrange multipliers.
So the structures of \reff{sosrlx1}-\reff{sosrlx2}
should be exploited to improve the efficiency of the method.
Section~5 discusses the specifications
with the nonnegative orthant $\re^n_+$ and
the standard box $[a, b]_n$.

{\bf Acknowledgments } \quad
The authors would like very much to thank Prof. Sturmfels for
the illuminating discussions with him and
his many constructive comments on the paper.


\end{document}